# ON THE SPACES OF POLYNOMIAL KNOTS

VICTOR A. VASSILIEV*
*Mathematics College*
*Independent Moscow University*
*Chernjakhovsky str. 5-1 apt. 150*
*Moscow 125319 Russia*

**Abstract**

Homology groups of spaces of nonsingular polynomial embeddings $\mathbf{R}^1 \to \mathbf{R}^n$ of degrees $\leq 4$ are calculated. A general algebraic technique of such calculations for spaces of polynomial knots of arbitrary degrees is described.

## 1. Introduction

The space of knots in $S^n$ with fixed behavior at the distinguished point can be approximated by the spaces of polynomial embeddings $\Phi : \mathbf{R}^1 \to \mathbf{R}^n$ of the form

$$\begin{cases} x_1(t) = t^d + a_1^1 t^{d-1} + \cdots + a_{d-1}^1 t \\ \cdots\cdots\cdots\cdots\cdots\cdots\cdots\cdots \\ x_n(t) = t^d + a_1^n t^{d-1} + \cdots + a_{d-1}^n t. \end{cases} \tag{1.1}$$

For instance, any $C^1$-knot in $S^3$ is isotopy equivalent to the closure of the image of such a map for some $d$. In [1] A. Shastri showed that the trefoil knot can be realized by such an embedding with $d = 5$. We prove in particular, that nontrivial knots cannot be realized by the maps of degrees less than 5.

Denote by $\mathcal{K}_d^n$ the space of all maps of the form (1.1), and by $\Sigma \equiv \Sigma_d$ *the discriminant* of $\mathcal{K}_d^n$, i.e. the closure of the space of maps (1.1) having either multiple or singular (with $\Phi'(t) = 0$) points. Our main result here is the following theorem.

**Theorem 1.** *For any $n \geq 2$, the space $\mathcal{K}_3^n \setminus \Sigma$ is contractible, and the space $\mathcal{K}_4^n \setminus \Sigma$ is homology equivalent to $S^{n-2}$. For $n \geq 4$ the space $\mathcal{K}_4^n \setminus \Sigma$ is also homotopy equivalent to $S^{n-2}$. The generator of the group $H^{n-2}(\mathcal{K}_4^n \setminus \Sigma)$ is equal to the linking number in $\mathcal{K}_4^n$ with the set of maps (1.1) having singular points.*

---

*Supported by the International Science Foundation (grant MQO000) and the Russian Fund of Basic Investigations. e.mail merx@glas.apc.org (For V.A.Vassiliev)

The first assertion of the theorem is trivial and should be known to the specialists. Indeed, it is easy to verify that the space $\mathcal{K}_3^n \setminus \Sigma$ can be contracted to the point (1.1) with $x_1(t) \equiv \cdots \equiv x_n(t) \equiv t^3 + t$ : if a map (1.1) belongs to $\mathcal{K}_3^n \setminus \Sigma$, then also the entire segment connecting it with this distinguished map does.

**Proposition 1.** *If $d$ is odd, then the space of maps (1.1) with multiple or singular points is closed and thus coincides with $\Sigma_d$. If $d$ is even, then $\Sigma_d$ is the union of this space and the space of all maps (1.1) with $a_1^1 = \cdots = a_1^n$.*

(For any $\alpha \in \mathbf{R}$ the space of maps (1.1) with $a_1^1 = \cdots = a_1^n = \alpha$ is the limit position of the spaces of maps (1.1) gluing together the points $t$ and $\alpha - t \in \mathbf{R}^1$, $t \to \infty$). □

**Proposition 2.** *If $n \geq 5$ or $n = d = 4$, then the space $\mathcal{K}_d^n \setminus \Sigma$ is simply-connected. (In particular, the homotopical part of Theorem 1 follows from the homological one by the Whitehead theorem.)*

Indeed, $\Sigma$ is a subvariety of codimension $n-2$, hence the proposition is true if $n \geq 5$, and in the case $n = 4$ the group $\pi_1(\mathcal{K}_d^n \setminus \Sigma)$ is generated by simple loops embracing $\Sigma$ at several nonsingular points of it. Moreover, as we shall see in § 4, if $d \geq 4$, then $\Sigma$ has a singularity-resolution $F_1\sigma \to \Sigma$, where $F_1\sigma$ is a connected smooth manifold with nonempty boundary, and the codimension in it of the pre-image of the singular set of $\Sigma$ is equal to $n-2$. Therefore any simple loop in $\mathcal{K}_d^n \setminus \Sigma$ embracing $\Sigma$ at nonsingular point $\xi \in \Sigma$ can be unhooked from $\Sigma$ by a deformation which follows a path connecting $\xi$ with $\partial\Sigma$ in the nonsingular part of $\Sigma$. □

**Proposition 3.** *For any even $d$,*

*a) the space $\mathcal{K}_d^3 \setminus \Sigma$ is a direct multiple of a circle;*

*b) for any $n \geq 3$, the group $H^{n-2}(\mathcal{K}_d^n \setminus \Sigma)$ contains at least one free generator, the linking number with the space of maps (1.1) with $a_1^1 = \cdots = a_1^n$;*

*c) this linking number is equal to the linking number with the (suitably oriented) set of maps (1.1) having a singular point.*

Indeed, the group $S^1$ of rotations of the space $\mathbf{R}^3$ around the axis $(1,1,1)$ acts freely on the space of maps (1.1) whose coefficients $a_1^1, a_1^2, a_1^3$ are not all equal, and any orbit of this action contains a unique map (1.1) with $a_1^1 = a_1^2$ and $a_1^1 + a_1^2 - 2a_1^3 > 0$. This proves assertion a). All these orbits have the linking number 1 with the set of maps with $a_1^1 = a_1^2 = a_1^3$. (For $d = 4$, the identical imbedding of any of these orbits into $\mathcal{K}_d^3 \setminus \Sigma$ induces the homology equivalence promised in Theorem 1.)

For arbitrary $n$, the linking number described in assertion b) is non-trivial because the orbit of the map (1.1) with $x_1(t) \equiv \cdots \equiv x_{n-1}(t) \equiv t^d$, $x_n(t) = t^d + t^{d-1} + t$ under the similar action of the group $SO(n-1)$ of rotations preserving the line $(1, \ldots, 1)$ is diffeomorphic to $S^{n-2}$ and the linking number takes value $\pm 1$ on it. Assertion c) will be proved in § 4.1.

Our method of calculation can be, in principle, applied to the cohomology of spaces of polynomial knots of arbitrary degree $d$. For $n = 3$ and $d \to \infty$, these cohomologies stabilize to the cohomology group of the space of knots, which includes the space of all knot invariants. The cocycles (in particular invariants) which stabilize well were introduced in [2], [3]. Here we start the regular investigation of arbitrary (nonstable) cocycles and the process of their stabilization.

# 2. The Conical Resolutions of Discriminants

As in the previous papers [2], [3], we exploit the idea of [4]: we consider the *closed homology group*

$$\bar{H}_*(\Sigma), \tag{2.1}$$

i.e. the homology group of the one-point compactification of $\Sigma$ reduced modulo the added point. The cohomology group of $\mathcal{K}_d^n \setminus \Sigma$ is Alexander dual to this group.

To calculate the group (2.1) we construct the *conical resolutions* of the discriminant varieties (used previously in [5]), which are a continuous generalization of the *simplicial resolutions* widely used in the algebraic geometry, see also [3], [6].

## 2.1. Primary conical resolution

In our situation, the simplest conical resolutions are constructed as follows. Denote by $P_d$ the space of polynomials of the form

$$t^d + a_1 t^{d-1} + \cdots + a_{d-1} t, \tag{2.2}$$

so that $\mathcal{K}_d^n \simeq (P_d)^n$. Consider the disjoint union of all Grassmann manifolds $G_i(P_d)$ of $i$-dimensional *affine* subspaces in $P_d$, $i = 0, \ldots, d-2$. The desired conical resolution is a subset of the direct product of $\mathcal{K}_d^n$ and the join of these Grassmannians,

$$G_0(P_d) * G_1(P_d) * \cdots * G_{d-2}(P_d). \tag{2.3}$$

**Definition 1.** An *elementary condition* on the polynomials $f \in P_d$ is either of the following three conditions:
a) for some distinguished points $t \neq s \in \mathbf{R}^1$, $f(t) = f(s)$;
b) for some distinguished point $t$, $f'(t) = 0$;
c) (only if $d$ is even) the coefficient $a_1$ is equal to some distinguished value $\alpha$.

The space of all elementary conditions is naturally topologized: for odd $d$, it is homeomorphic to the closed half-plane

$$\{(t, s) \in \mathbf{R}^2 | t \leq s\}, \tag{2.4}$$

where the diagonal $\{t = s\}$ corresponds to the conditions b). For even $d$, this space is homeomorphic to the infinite strip $\mathbf{R}^1 \times [0, 1]$ and is obtained from the above half-plane by adding a line "at infinity": the point $\{\alpha\}$, $\alpha \in \mathbf{R}$, of this line is the common limit point of all non-compact curves in the half-plane (2.4) tending to infinity in such a way that the values of the function $t + s$ on them converge to $\alpha$ when $t$ tends to $-\infty$ and $s$ to $+\infty$. This line is called the line of *exceptional conditions.*

In both cases, we denote this space by $\Psi(d)$.

If $d > 2$, then any elementary condition defines a point in $G_{d-2}(P_d)$: the space of all polynomials satisfying it. Any point in $G_{d-2}(P_d)$ corresponds to at most one elementary condition; the above-described topology in the space $\Psi(d)$ is no other but that induced by this correspondence from the standard topology in $G_{d-2}(P_d)$.

Any collection $T$ of elementary conditions defines an affine subspace $L(T)$ in $P_d$, the set of polynomials satisfying all these conditions, and thus a point in $G_i(P_d)$ for some $i$.

Denote by $C(d, i)$ the set of all $i$-dimensional affine subspaces of the form $L(T)$ in $P_d$, let $\bar{C}(d, i)$ be the closure of it in $G_i(P_d)$.

The join (2.3) is, roughly speaking, the naturally topologized union of all simplices whose vertices belong to manifolds $G_i(P_d)$ with different $i$.

**Definition 2.** A *coherent simplex* in this join is a simplex (of arbitrary dimension) whose vertices are the points of subsets $\bar{C}(d, i) \subset G_i(P_d)$ with different $i$, and **all** the planes in $P_d$ corresponding to these vertices are incident to one another.

Denote by $\Lambda(d)$ the union of all non-empty coherent simplices in the join (2.3).

To any simplex $\Delta$ constituting $\Lambda(d)$ there corresponds an affine plane $L(\Delta) \subset P_d$, namely, the vertex of $\Delta$ lying in $G_i(P_d)$ with the smallest value of $i$, and hence also a plane $\mathcal{L}(\Delta)$ in $\mathcal{K}_d^n$: the set of maps (1.1) all $n$ components of which belong to $L(\Delta)$.

The conical resolution of the discriminant, $R\Sigma$, is defined as the subset in $\Lambda(d) \times \mathcal{K}_d^n$, which is the union of pairs of the form (an interior point of some coherent simplex $\Delta$, a point of the corresponding plane $\mathcal{L}(\Delta) \subset \mathcal{K}_d^n$).

Define the map $\pi : R\Sigma \to \Sigma$ as the restriction on $R\Sigma$ of the obvious projection $\Lambda(d) \times \mathcal{K}_d^n \to \mathcal{K}_d^n$.

**Proposition 4.** *For any $d$, the map $\pi$ is proper and surjective. It defines a homotopy equivalence of one-point compactifications of spaces $R\Sigma$ and $\Sigma$. $R\Sigma$ is homeomorphic to a semialgebraic set in such a way that $\pi$ becomes a piecewise-algebraic map.*

Indeed, the first assertion is obvious, and the last one follows from the construction if we realize the join (2.3) as the union of simplices spanning the points of Grassmannians $G_i(P_d)$ imbedded algebraically and generically into the Euclidean space of a very large dimension. The second assertion follows now from the fact that the pre-image of any point $\phi \in \Sigma$ is a contractible compact $CW$-complex, namely, the union of all coherent simplices $\Delta$ whose common principal vertex $L(\Delta)$ is the minimal plane in $P_d$ belonging to some $\bar{C}(d, i)$ and containing simultaneously all $n$ components of the map $\phi$. $\Box$

## 2.2. Reduced conical resolution

Usually the resolution constructed above is not the most convenient one and can be reduced a bit to avoid the unnecessary calculations.

**Example.** A generic plane $L(T)$ of codimension two in $P_d$, $d \geq 4$, is distinguished by two different elementary conditions $x, y \in \Psi(d)$ of type a) from Definition 1. Such a point $L(T) \in G_{d-3}(P_d)$ is connected by coherent segments with exactly two planes $L(x), L(y) \in G_{d-2}(P_d)$. However, when $x$ tends to $y$, the limit value of $L(x,y)$ is connected with unique point $L(y) \in G_{d-2}(P_d)$; this limit value belongs to $\bar{C}(d, d-3)$ but not to $C(d, d-3)$. (Such limit values for a given $y$ run over a circle and depend on the angle of attack of $x$ towards $y$; these limit planes in $P_d$ sweep out the plane $L(y)$ and define a blowing-up of it with the center at the subspace of polynomials with vanishing derivatives at both points constituting the pair $y$).

In general, for any plane $L \in \bar{C}(d,i)$ denote by $\bar{L}$ the minimal plane of the form $L(T)$ containing $L$: by definition, this is a point of $C(d,j)$ for some $j \geq i$. Conversely, to any point $L' \in C(d,j)$ there corresponds a subcomplex in $\Lambda(d)$ consisting of all points $L \in \bar{C}(d,i)$, $i \leq j$, such that $L' = \bar{L}$, and of all simplices spanning such points with different $i$. Obviously this subcomplex is compact and contractible.

The *reduction map red* : $\Lambda(d) \to \tilde{\Lambda}(d)$ is defined as follows: for any $i$ and any $L \in \bar{C}(d,i)$ we map the point $L$ to $\bar{L}$ and extend this map by linearity to all coherent simplices; the image of a coherent simplex under this map is obviously again coherent. Then we identify to one another all the points of $\Lambda(d)$ sent by this extended map to the same point. $\tilde{\Lambda}(d)$ is the quotient space of $\Lambda(d)$ with respect to this identification.

This map defines in obvious way a map $\Lambda(d) \times \mathcal{K}_d^n \to \tilde{\Lambda}(d) \times \mathcal{K}_d^n$ and hence also a map of the resolved discriminant $R\Sigma$ onto some subcomplex $\tilde{R}\Sigma \subset \tilde{\Lambda}(d) \times \mathcal{K}_d^n$.

**Proposition 5.** *The extended identification map*

$$Red : R\Sigma \to \tilde{R}\Sigma \qquad (2.5)$$

*induces the homotopy equivalence of the one-point compactifications of $R\Sigma$ and $\tilde{R}\Sigma$ and factorizes the map $\pi : R\Sigma \to \Sigma$ (i.e. there is a proper map $\tilde{\pi} : \tilde{R}\Sigma \to \Sigma$ such that $\pi = \tilde{\pi} \circ Red$).*

This follows immediately from the construction. □

## 2.3. Filtration and spectral sequence

The space $\Lambda(d)$ admits a natural increasing filtration. Indeed, the filtration of any interior point of a simplex $\Delta \subset \Lambda(d)$ is equal to the codimension $d-1-i$ of the principal plane $L(\Delta) \in \bar{C}(d,i)$ of $\Delta$. Hence also a filtration on the space $\tilde{\Lambda}(d)$ is induced: the filtration of a coset $x \in \tilde{\Lambda}(d)$ is equal to the smallest value of filtrations of points in $\Lambda(d)$ constituting $x$. The obvious projections $R\Sigma \to \Lambda(d)$, $\tilde{R}\Sigma \to \tilde{\Lambda}(d)$ induce the filtrations also on the spaces $R\Sigma$ and $\tilde{R}\Sigma$. By the construction, the term $F_i \setminus F_{i-1}$ of this filtration of $R\Sigma$ (respectively, $\tilde{R}\Sigma$) is the space of a $n(d-1-i)$-dimensional affine bundle over the similar term of the filtration of $\Lambda(d)$ (respectively, $\tilde{\Lambda}(d)$).

This filtration defines a spectral sequence converging to the group (2.1): the term $E^1_{p,q}$ of it is isomorphic to $\bar{H}_{p+q}(F_p \setminus F_{p-1})$. For some reasons (see f.i. [3]) it is convenient to rename this term to $E_1^{-p,n(d-1)-1-q}$ so that the spectral sequence becomes a cohomological one; this reversed spectral sequence converges exactly to the cohomology of $\mathcal{K}^n_d \setminus \Sigma$.

# 3. The Trivial Example

Although the topology of the space $\mathcal{K}^n_3 \setminus \Sigma$ is clear, we now calculate its cohomology to demonstrate the general method.

**Definition 3.** An *open cone* over a topological space $X$ is the usual cone $CX$ with its base $X$ deleted.

The term $F_1$ of the natural filtration of $\tilde{\Lambda}(3)$ (and of $\Lambda(3)$) is the space $C(3,1) \equiv \Psi(3)$, and the similar term of the filtration of $\tilde{R}\Sigma$ is the space of a $n$-dimensional affine bundle over $\Psi(3)$. (Namely, any element $x \in \Psi(3)$ defines a (one-dimensional) pencil of polynomials in $P_3$; the fibre over $x$ of the promised bundle is the direct sum of $n$ examples of this pencil.) In particular, $E_1^{-1,*} \equiv 0$.

The spaces $L(T)$ of codimension two in $P_3$ are just the points, thus to calculate the term $E_1^{-2,*}$ we consider all the polynomials (2.2) of degree $d = 3$ and for any such polynomial $f$ take a cone whose base is the space of all affine lines $L(x) \in C(3,1) \equiv \Psi(3)$ containing $f$; the term $F_2 \setminus F_1$ of $R\Sigma$ is swept out by all such (open) cones (the bases of them lie in the term $F_1$). To do this we need the following notion.

**Definition 4.** The *relation curve* $r(f)$ defined by a polynomial $f : \mathbf{R}^1 \to \mathbf{R}^1$ is the set of all pairs $(t,s) \in \Psi(d)$ such that $f(t) = f(s)$ (if $t \neq s$) or $f'(t) = 0$ (if $t = s$); if $d$ is even, then this curve includes also the point of the exceptional line in $\Psi(d)$ equal to the coefficient $a_1$ of $f$.

**Lemma 1.** *If $f \in P_d$, $d$ odd, then the curve $r(f)$ lies in the "finite" domain in $\Psi(d)$. If $d$ is even, then $r(f)$ has unique branch going to the infinity, and the line $\{t + s = a_1\}$ is asymptotical for it.* □

The polynomials $f_1, \ldots, f_k$ belong to a subspace $L(T)$ (see § 2.1) if and only if the set $r(f_1) \cap \ldots \cap r(f_k)$ contains all elementary conditions $w \in T$.

The relation curve of a cubic polynomial $f = t^3 + \alpha t^2 + \beta t$ is distinguished by the equation $(t^2 + ts + s^2) + \alpha(t+s) + \beta = 0$, hence is either empty (if $f'$ is everywhere positive), consists of unique point (if $f'$ has a double root) or is a half-ellipse in the half-plane $\Psi(3)$ orthogonal to the boundary of $\Psi(3)$ at its two endpoints (which are such diagonal points $\{t,t\}$ that $t$ is a root of $f'$).

The polynomials with $f' > 0$ (which are distinguished by the inequality $\alpha^2 < 3\beta$) obviously do not enter any pencil sweeping $F_1$.

Let $f$ be any polynomial with two real critical points. Then for any point of the corresponding relation curve $r(f)$ there is exactly one pencil satisfying it. Hence the subset in $F_2 \setminus F_1$ corresponding to such a polynomial $f$ is an open cone over a half-ellipse ($\sim$ a segment): the segment in the base of this cone belongs to $F_1$. Thus the part of $F_2 \setminus F_1$ swept out by the cones corresponding to such $f$ is a fibre bundle over the 2-cell $\{(\alpha,\beta) \in \mathbf{R}^2 | \alpha^2 > 3\beta\}$ with fibre homeomorphic to a triangle with one closed side deleted. The (closed) homology of this part is trivial.

Finally, for border values of $f$ (i.e. $\alpha^2 = 3\beta$) the base segment of such open cones contracts to a point, thus the cone corresponding to such an $f$ in the term $F_2 \setminus F_1$ of the natural filtration of $R\Sigma$ is a half-open interval; after going to $\tilde{R}\Sigma$ these intervals vanish and give trivial deposit to the homology of the discriminant. Thus the column $E_1^{-2,*}$ also vanishes and the whole group $\bar{H}_*(\Sigma)$ is trivial.

# 4. Space of Polynomial Knots of Fourth Degree

## 4.1. The column $E_1^{-1,*}$

For any $d$, the term $F_1$ of the standard filtration of $\tilde{R}\Sigma$ is the space a $n(d-2)$-dimensional affine bundle over $\Psi(d)$, thus for even $d$ the group $E_1^{-1,q}$ is free cyclic if $q = n-1$ and is trivial for all other $q$. The generator of the group $E_1^{-1,n-1}$ can be realized as the fundamental class of the (unfolded) set of maps (1.1) with singular points or of the (homological to it) set of maps with $a_1^1 = \cdots = a_1^n$.

This gives also a proof of assertion c) of Proposition 3. □

## 4.2. The column $E_1^{-2,*}$

The main result of this subsection is the following theorem.

**Theorem 2.** $E_1^{-2,*} \equiv 0$.

By the general algorithm, we should consider all one-dimensional affine subspaces in $P_4$ and for any such line take an open cone over the set of all elementary conditions $x \in \Psi(4)$ satisfied by all polynomials of this line. The term $F_2 \setminus F_1$ of $R\Sigma$ is the space of a $n$-dimensional affine bundle over the space swept out by all these cones.

**Lemma 2.** *The number of elementary conditions satisfied by all polynomials of some affine line in $P_d$ can be equal to* $0, 1, \ldots, (d-1)(d-2)/2$ *or* $\infty$.

Let us represent these lines in canonical form. A line is determined by the choice of any its point $f$ and any directing vector $g$, i.e. the difference of arbitrary two distinct points. This direction $g$ is a polynomial of degree strictly lower than $d$.

**Definition 5.** The pair $(f, g)$ representing a line in $P_d$ is chosen *in canonical form* if the leading coefficient of $g$ is equal to 1 and the coefficient of the monomial $t^{\deg\, g}$ in $f$ is equal to 0.

Obviously, any affine line in $P_d$ can be uniquely represented by a pair of polynomials $(f, g)$ in canonical form. The set of elementary conditions satisfied by all the polynomials of the line is nothing but the set $r(f) \cap r(g)$.

The relation curve of a polynomial of degree $d$ has degree $d-1$, hence Lemma 2 follows from the Bezout's theorem and the fact that the relation curves are invariant under the reflection $(t, s) \to (s, t)$ and we count only the (geometrical) points with $t \leq s$. $\Box$

In particular, for $d = 4$, the cardinality of $r(f) \cap r(g)$ can be equal to $0, 1, 2, 3$ or $\infty$. Obviously, if this number is infinite, then the curve $r(f)$ is reducible.

**Lemma 3.** *For $f$ of degree 4, the curve $r(f)$ is reducible only if the polynomial $f$ is symmetric with respect to some point $t_0 \in \mathbf{R}^1$, i.e., $f = (t-t_0)^4 + \lambda(t-t_0)^2 + c$ for some $\lambda$ and $c = -\lambda t_0^2$. In this case $r(t)$ consists of the segment $\{(t, s)|t+s = 2t_0\}$ and the half-circle (maybe imaginary) $(t-t_0)^2 + (s-t_0)^2 + \lambda = 0$.* $\Box$

Obviously such a circle cannot coincide with a relation curve of a polynomial of degree 3, hence the pair $(f, g)$ satisfies infinitely many elementary conditions only if $\deg g = 2$ and the origin of symmetry of $g$ in $\mathbf{R}^1$ coincides with that of $f$.

If our line $(f, g)$ in $P_4$ satisfies 0 or 1 elementary conditions in $\Psi(4)$, then over it there are no points in the term $F_2 \setminus F_1$ of $\tilde{R}\Sigma$ (in the first case also in similar term of $R\Sigma$); if it satisfies 2 or 3 conditions, then the corresponding set in $F_2 \setminus F_1$ is an open interval (respectively, a three-ray star without ends).

### 4.2.1. Revised conical resolution

Now we slightly improve the construction of the conical resolution of $\Sigma \subset P_4$ replacing it (and any term $F_i \setminus F_{i-1}$ of its natural filtration) by a homology equivalent space. For instance, the three-end stars without endpoints which corresponded in the canonical construction of $\Lambda(4)$ to any line in $P_4$ satisfying exactly three conditions will be replaced by the triangles without vertices.

**Lemma 4.** *For any pair of distinct elementary conditions in $\Psi(4)$, which are not both exceptional, the set of polynomials in $P_4$ satisfying these two conditions is an affine line. All affine lines in $P_4$ whose points satisfy simultaneously at least two elementary conditions can be obtained in this way.* $\Box$

Denote by $\Psi^{[2]}(4)$ the space of unordered pairs of distinct elementary conditions, at least one of which is non-exceptional. So we have a continuous map

$$\nu : \Psi^{[2]}(4) \to \bar{C}(4, 1) \qquad (4.1)$$

that sends a pair of conditions $(x, y)$ to the pencil of polynomials satisfying them. Let $\bar{\Psi}^{[2]}(4)$ be the completion of $\Psi^{[2]}(4)$ induced by this map from the completion $\bar{C}(4,1)$ of $C(4,1)$ (i.e. $\bar{\Psi}^{[2]}(4)$ is the union of $\Psi^{[2]}(4)$ and the set $\partial\Psi^{[2]}(4)$ canonically homeomorphic to $\bar{C}(4,1) \setminus C(4,1)$, topologized in such way that a sequence $\{x_i\} \in \Psi^{[2]}(4)$ converges to a point $y \in \partial\Psi^{[2]}(4)$ if and only if the sequence $\{\nu(x_i)\} \in C(4,1)$ converges to the point of $\bar{C}(4,1) \setminus C(4,1)$ corresponding to $y$ via the canonical homeomorphism). The obvious extended map $\bar{\Psi}^{[2]}(4) \to \bar{C}(4,1)$ will be also denoted by $\nu$. Now we consider the join

$$\bar{C}(4,0) * \bar{C}(4,1) * \bar{\Psi}^{[2]}(4) * \bar{C}(4,2). \tag{4.2}$$

To any point of $\bar{C}(4,0)$, $\bar{C}(4,1)$, $\bar{\Psi}^{[2]}(4)$ or $\bar{C}(4,2)$ there corresponds an affine subspace in $P_4$: for spaces $\bar{C}(4,i)$ this correspondence is tautological, and to the point $x \in \bar{\Psi}^{[2]}(4)$ there corresponds the line $\nu(x)$. A simplex of arbitrary dimension in the join (4.2) is called *coherent*, if all planes in $P_4$ corresponding to its vertices are incident to one another (in particular, if this simplex contains the points of both $\bar{C}(4,1)$ and $\bar{\Psi}^{[2]}(4)$, then the corresponding lines should coincide).

The forthcoming construction of the conical resolutions of $\Sigma$ based on the union of coherent simplices in (4.2) repeats that from §§ 2.1, 2.2; the spaces of these resolutions will be denoted respectively by $\rho\Sigma$ and $\tilde{\rho}\Sigma$, the union of coherent simplices by $\lambda(4)$ and its natural quotient space (cf. § 2.2) by $\tilde{\lambda}(4)$.

The join (4.2) (and hence also its subspace $\lambda(4)$) admits the following increasing filtration: $F_1$ consists of the points of $C(4,2)$, $F_{1.5}$ consists of points of $C(4,2)$ and $\bar{\Psi}^{[2]}(4)$ and segments connecting them; $F_2$ consists of the points of $\bar{C}(4,1)$, $\bar{\Psi}^{[2]}(4)$ and $C(4,2)$ and the segments and triangles spanning them; finally, $F_3$ is the entire join (4.2). As in § 2.3, this filtration induces filtrations on $\tilde{\lambda}(4)$, $\rho\Sigma$ and $\tilde{\rho}\Sigma$.

There are natural filtration-preserving projections

$$\lambda(4) \to \Lambda(4), \quad \tilde{\lambda}(4) \to \tilde{\Lambda}(4). \tag{4.3}$$

For instance, the first of them maps linearly any coherent simplex in (4.2) on its face belonging to $\bar{C}(4,0) * \bar{C}(4,1) * \bar{C}(4,2)$, is identical on this face and maps the vertex $x \in \bar{\Psi}^{[2]}(4)$ to the vertex $\nu(x) \in \bar{C}(4,1)$. The second projection in (4.3) is induced by the first one. These projections define also the maps

$$\rho\Sigma \to R\Sigma, \quad \tilde{\rho}\Sigma \to \tilde{R}\Sigma; \tag{4.4}$$

indeed, the maps (4.3) can be lifted in obvious way to the maps $\lambda(4) \times \mathcal{K}_4^n \to \Lambda(4) \times \mathcal{K}_4^n$, $\tilde{\lambda}(4) \times \mathcal{K}_4^n \to \tilde{\Lambda}(4) \times \mathcal{K}_4^n$; the desired maps (4.4) are the restrictions of them on the subsets $\rho\Sigma \subset \lambda(4) \times \mathcal{K}_4^n$, $\tilde{\rho}\Sigma \subset \tilde{\lambda}(4) \times \mathcal{K}_4^n$.

**Proposition 6.** *The projections (4.4) define the homotopy equivalences of the one-point compactifications of all involved spaces. Moreover, for any $i = 1,2,3$ the same is true for the induced maps of the terms $F_i \setminus F_{i-1}$ of our filtrations of the spaces $\rho\Sigma$ (respectively, $\tilde{\rho}\Sigma$) to $R\Sigma$ (respectively, $\tilde{R}\Sigma$). In particular, these maps induce isomorphisms of the corresponding spectral sequences.* □

(The additional term $F_{1.5}$ of the filtrations of $\rho\Sigma$ and $\tilde{\rho}\Sigma$ is just a mean in calculating the groups $E_1^{-2,i}$ of them.) Now Theorem 2 follows from the next two lemmas (concerning the above-described filtration of the space $\tilde{\rho}\Sigma$).

**Lemma 5.** $\bar{H}_*(F_{1.5} \setminus F_1) \equiv 0$.

**Lemma 6.** $\bar{H}_*(F_2 \setminus F_{1.5}) \equiv 0$.

*Proof of Lemma 5.* The term $F_{1.5} \setminus F_1$ of the natural filtration of $\tilde{\lambda}(4)$ is homeomorphic to a fibre bundle over $\Psi^{[2]}(4)$ with an interval as a fibre: indeed, to any point $\{x, y\} \in \Psi^{[2]}(4)$ there corresponds the union of two half-open intervals $[\{x, y\}; x)$ and $[\{x, y\}; y)$. It is easy to calculate that the unique nontrivial element of the fundamental group of $\Psi^{[2]}(4)$ changes the orientations of this bundle of intervals and of any of $n$ line factors of the canonical $n$-dimensional affine bundle $(F_{1.5} \setminus F_1)(\tilde{\rho}\Sigma) \to (F_{1.5} \setminus F_1)(\tilde{\lambda}(4))$. Thus the group $\bar{H}_i((F_{1.5} \setminus F_1)(\tilde{\rho}\Sigma))$ is isomorphic to $\bar{H}_{i-n-1}(\Psi^{[2]}(4))$ if $n$ is even and to the group $\bar{H}_{i-n-1}(\Psi^{[2]}(4); \pm\mathbf{Z})$ (where $\pm\mathbf{Z}$ is the unique non-trivial local system with fibre $\mathbf{Z}$) if $n$ is odd. Both these groups can be easily calculated and are trivial. □

## 4.3. Proof of Lemma 6.

The term $F_2 \setminus F_{1.5}$ of $\tilde{\lambda}(4)$ consists of two parts "3" and "$\infty$" swept out by coherent triangles

$$[(f, g) \in G_1(P_4); \{x, y\} \in \Psi^{[2]}(4); w \in \Psi(4)] \tag{4.5}$$

such that the corresponding line $(f, g)$ satisfies exactly 3 (respectively, infinitely many) elementary conditions. The pre-images of these two sets under the obvious projection $\tilde{\rho}\Sigma \to \tilde{\lambda}(4)$ will be denoted by [3] and [$\infty$]. Now we calculate separately the closed homology of any of these two parts.

### 4.3.1. Cohomology of [3].

To any line $(f, g)$ satisfying exactly three conditions there correspond six triangles (4.5) in $\tilde{\lambda}(4)$, the union of which is an open triangle with barycenter at the point $(f, g)$: the boundary of this triangle belongs to $F_{1.5}$ (and the vertices even to $F_1$).

**Proposition 7.** *The set of lines $(f, g) \in G_1(P_4)$, all elements of which satisfy exactly three common elementary conditions, is homeomorphic to the direct product of the line $\mathbf{R}^1$, the interval $(-\infty, 0)$ and the plane set shown in Fig. 1a) (where the dashed segments denote the part of the boundary which do not belong to the set). In particular the closed homology of this set is trivial.*

**Corollary.** $\bar{H}_*([3]) = 0$.

Indeed, [3] is the space of a fibre bundle over the set in $G_1(P_4)$ considered in the previous proposition, the fibre of which is a direct product of a $n$-dimensional affine space (the $n$-th power of the tautological bundle) and an open triangle. □

*Proof of Proposition 7.* By the Bezout's theorem, if the pair $(f, g)$ satisfies exactly three conditions, then deg $g = 3$, $g = t^3 + \alpha t^2 + \beta t$ and $r(g)$ has the form

$$(t^2 + ts + s^2) + \alpha(t + s) + \beta = 0. \tag{4.6}$$

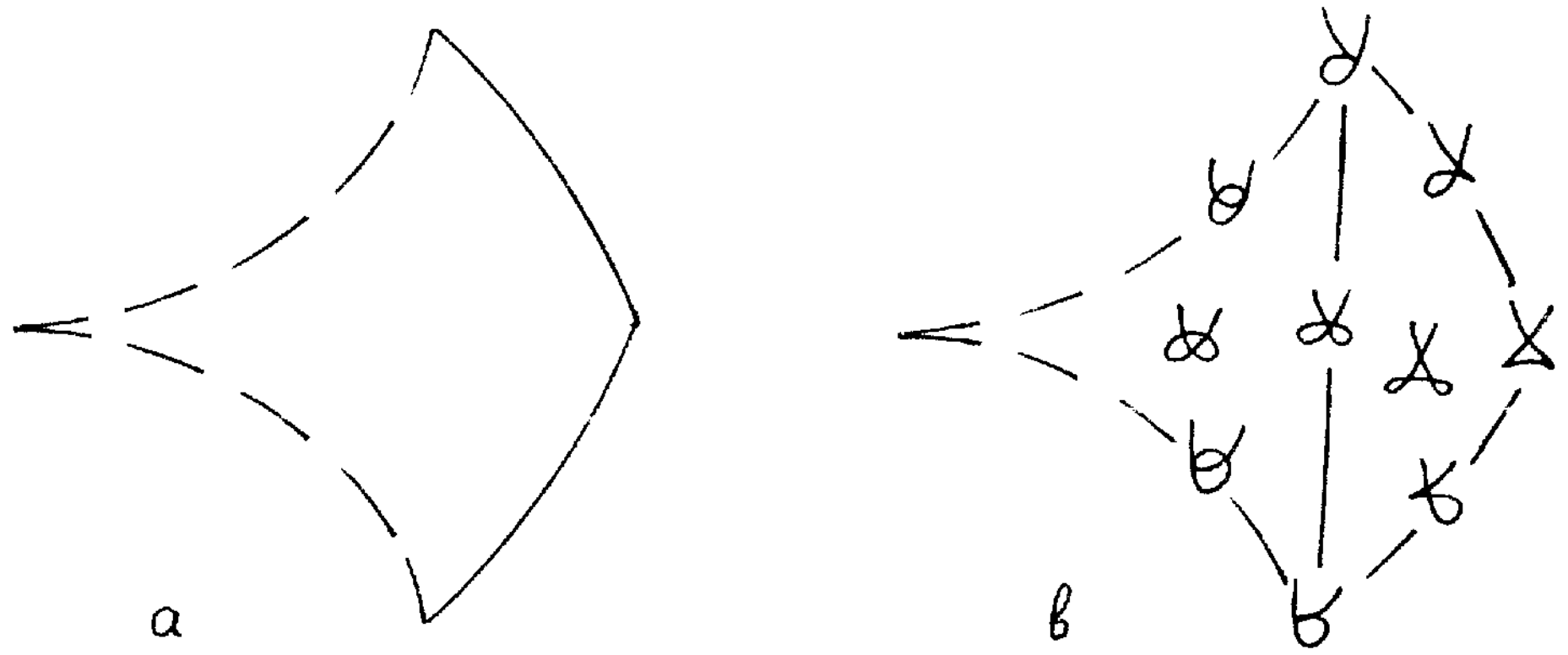

Fig. 1: The set of polynomials $f$ whose relation curve intersects $r(g)$ at three points

The group of translations in the line of arguments $t$ acts freely on the investigated set: this action defines the first factor $\mathbf{R}^1$ in Proposition 7. Hence it is sufficient to consider only the space of pairs $(f,g)$ with the coefficient $\alpha$ in $g$ equal to zero: this space can be naturally identified with the set of orbits of this action. The coefficient $\beta$ in (4.6) should obviously be negative; the values of it define the factor $(-\infty,0)$ in Proposition 7. Finally, we need only to prove the following lemma.

**Lemma 7.** *For any negative value $\beta$ the set of polynomials $f = t^4+at^2+bt$ such that the corresponding relation curve*

$$r(f) \equiv \{(t,s)|(t^3+t^2s+ts^2+s^3)+a(t+s)+b=0\} \qquad (4.7)$$

*intersects the curve*

$$r(g) \equiv \{(t,s)|(t^2+ts+s^2)+\beta=0\} \qquad (4.8)$$

*at exactly three points is homeomorphic to the set shown in Fig. 1a). The union of such sets over different values of $\beta$ defines a locally trivial fibre bundle over the space $(-\infty,0)$ of all such values.*

*Proof.* The group $\mathbf{R}_+$ acts on the space of pairs $(f,g)$ by the rule

$$T_\lambda(f(t),g(t)) = (\lambda^{-4}f(\lambda t), \lambda^{-3}g(\lambda t)); \qquad (4.9)$$

this action reduces the proof of the first assertion of Lemma 7 to the case $\beta=-1$ and proves the second assertion about the local triviality.

The pair of equations (4.7), (4.8) with $\beta=-1$, being rewritten in symmetric functions $p=t+s$ and $q=ts$, becomes the system

$$\begin{cases} p(p^2-2q)+ap+b=0 \\ p^2-q=1 \\ p^2-4q\geq 0, \end{cases} \qquad (4.10)$$

where the inequality stays for the reality of solutions. Avoiding $q$, we get the system

$$p^3 - p(a+2) - b = 0 \tag{4.11}$$

$$p^2 \leq 4/3. \tag{4.12}$$

We are interested in the values of $a$ and $b$ for which this system has exactly three real solutions. Let $A = a + 2$, then the set of pairs $(A, b)$ for which the equation (4.11) has three real roots is the domain bounded by the semicubical parabola,

$$27b^2 \leq 4A^3, \tag{4.13}$$

see Fig. 1a). This part is swept out by the parametric curves $\{A = \tau^2,\ b = C \cdot \tau^3\}$, $\tau \in (0, +\infty)$ with $C < \sqrt{4/27}$. For any such curve, the squares of all three roots of (4.11) grow monotonically along it, therefore this curve leaves the set distinguished by (4.12) once and for all. The boundary of this set inside the domain (4.13) consists of two components depending on the sign of the root of (4.11) that leaves this domain at the corresponding point: these two components meet on the line $\{b = 0\}$ where the absolute values of two roots coincide.

This part of the boundary of the set of values $\{a, b\}$ distinguished by (4.11), (4.12) belongs to the domain studied in Lemma 7: indeed, its points correspond to the pairs of polynomials $(f, g)$ having one (or two in the breakpoint) common singular points and two (respectively, one) selfintersection point. (The topological shape of the plane curves given by the maps $(f, g) : \mathbf{R}^1 \to \mathbf{R}^2$ for different values of $A$ and $b$ is given in Fig. 1b).) On the other hand, the part of the boundary on which the inequality (4.13) becomes an equality does not belong to the domain we are interested in: at the points of it two selfintersection points coalesce and become one point of self-tangency. This proves Lemma 7.

**Corollary.** *If all polynomials of the affine line $(f, g) \in P_4$ satisfy exactly three common elementary conditions, then for any $\lambda \in \mathbf{R}$*

*1) the polynomial $f + \lambda g$ of this line has three real critical points (at least two of which are geometrically distinct). Any segment in $\mathbf{R}^1$ connecting two critical points of such polynomial $f + \lambda g$ contains a critical point of $g$;*

*2) the topological disposition of the relation curves $r(f + \lambda g), r(g)$ is as shown in Figs. 2a), b) or in pictures symmetric to them with respect to the reflection in the vertical line $\{t + s = const\}$, or is some degeneration of them, corresponding to either the symmetric polynomials $f + \lambda g$ (see Figs. 2 c), d)), or to the coincidence of some critical points of $f$ with these of $g$ (in the latter case $r(f + \lambda g)$ and $r(g)$ intersect at their endpoints).*

Indeed, for $\lambda = 0$ and $g$ of the form $t^3 + \beta t$ the assertion 1 follows from the path-connectedness of the set investigated in Lemma 7. For arbitrary $\lambda$ the same follows from the Roll theorem since at the roots of $g'$ the signs of all polynomials $(f + \lambda g)'$ are the same. Finally, the case of arbitrary $\lambda$ and $g$ can be reduced to the previous one by the suitable translation of the argument line. Statement 2 follows immediately from 1.

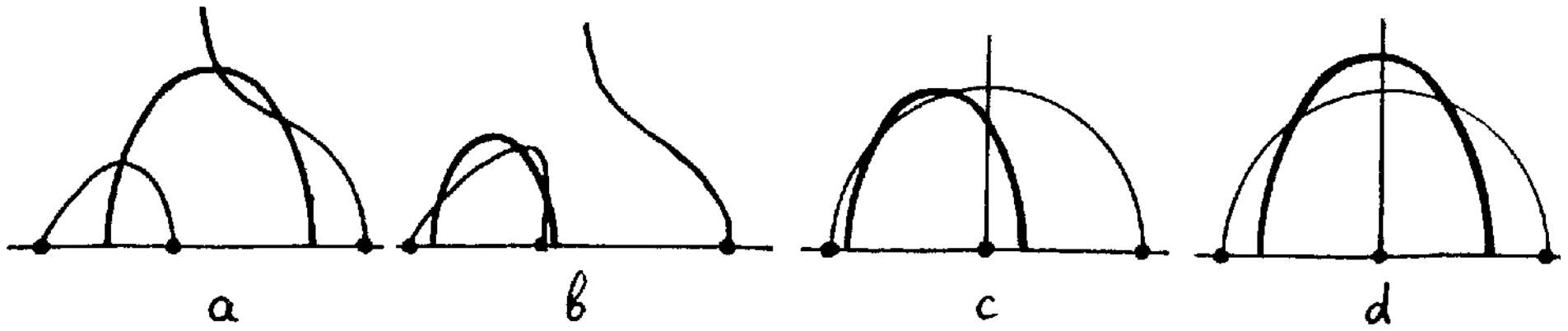

Fig. 2: Relation curves with three intersection points

### 4.3.2. Closed homology of $[\infty]$.

Consider a manifold $S$ and imbed it generically into the space $\mathbf{R}^N$ of a sufficiently large dimension in such a way that none two segments spanning the images of points of $S$ intersect in their interior points. Denote the union of such segments by $S^{*2}$. The space $S$ is naturally imbedded in it.

**Lemma 8.** *If $S$ is compact, then $S^{*2}$ is compact and the homotopy type of the pair $(S^{*2}, S)$ depends only on that of $S$.* □

**Proposition 8.** *The part "$\infty$" of $\tilde{\lambda}(4)$ is homeomorphic to the direct product of the line $\mathbf{R}^1$ and an open cone over the space $S^{*2}$ where $S$ is the segment $[0, +\infty]$.*

Indeed, if the set $r(f) \cap r(g)$ is infinite, then by Lemma 3 $\deg g = 2$ and this set is the compact segment in $\Psi(4)$ given by equation $t + s = 2t_0$ for some $t_0$; this coefficient $t_0$ defines the first factor $\mathbf{R}^1$ in Proposition 8. The structure of the cone follows from the construction of the resolution. □

**Corollary.** $\bar{H}_*([\infty]) = 0$.

Indeed, $[\infty]$ is the space of a $n$-dimensional bundle over "$\infty$", therefore $\bar{H}_i([\infty]) \equiv \bar{H}_{i-n}(\text{“}\infty'')$; by Proposition 8 the latter group coincides with the closed homology group of a direct product of a line and a half-open interval. □

Theorem 2 is completely proved. □

## 4.4. The column $E_1^{-3,*}$.

In this subsection we prove the following theorem.

**Theorem 3.** $E_1^{-3,*} \equiv 0$.

On this step, we consider all the subspaces of codimension three in $P_4$, i.e. the single points of it. For any such point $f$ we consider the pre-image in $F_2$ of the map $(f, \ldots, f) \in \mathcal{K}_4^n$ under the obvious projection $\pi : \rho\Sigma \to \Sigma$ and take the open cone over this pre-image. Such cones over all $f \in P_4$ sweep out the term $F_3 \setminus F_2$ of the resolution $\rho\Sigma$; similar term of $\tilde{\rho}\Sigma$ is swept out by such cones over those $f$ which belong to the subset $C(4, 0) \subset P_4$.

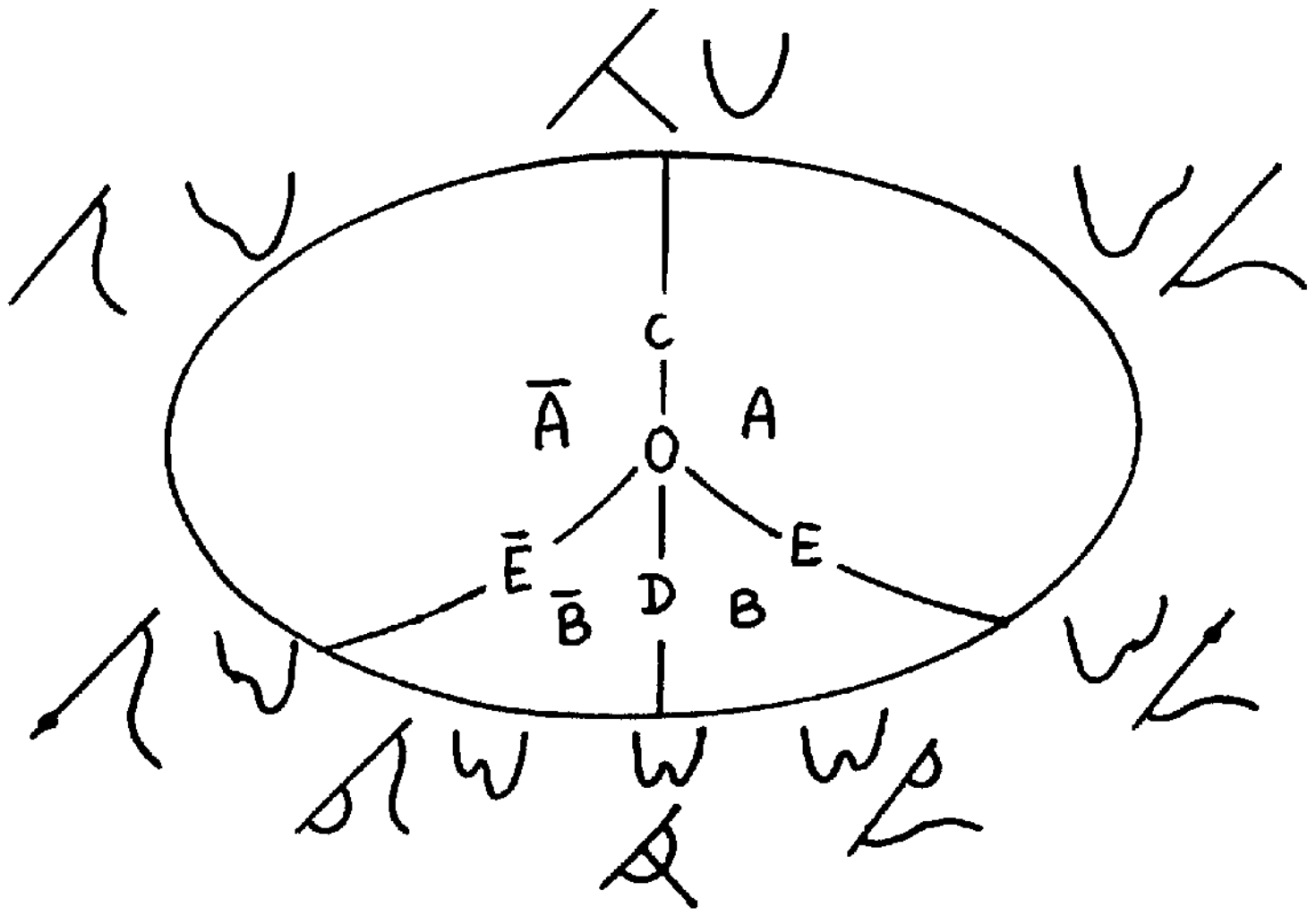

Fig. 3: Classification of polynomials of fourth degree

Using as before the free action of the group of translations of the argument space we can consider only the polynomials $f$ of the form $f = t^4 + at^2 + bt$, i.e. with vice-leading coefficient equal to 0: the entire space $F_3 \setminus F_2$ is the direct product of the line $\mathbf{R}^1$ and the subspace swept out by cones corresponding to such $f$.

The space of such polynomials is a two-plane and is naturally subdivided into 9 open cells, see Fig. 3. For any cell, we show in this picture the shape of the corresponding polynomials and relation curves. The cells $A$ and $A'$ (respectively, $B$ and $B'$) contain the strictly Morse polynomials with one (respectively, three) real critical points. They are separated by the semicubical parabola $\{27b^2 + 8a^3 = 0\} \equiv E \cup E' \cup O$ and the vertical line $\{b = 0\} \equiv D \cup C \cup O$ consisting of even polynomials.

For any of the nine cells, we consider the union of maps $(f, \dots, f) \in \Sigma \subset \mathcal{K}_4^n$ where $f$ runs over the cell, and consider the pre-image in $F_3 \setminus F_2$ of this union under the projection

$$\tilde{\pi} : \tilde{\rho}\Sigma \to \Sigma. \tag{4.14}$$

Denote this pre-image by $\{X\}$ where $X$ is the notation of the cell, $X = A, A', \dots$

### 4.4.1. The cells $C$ and $O$.

All points $f$ of the ray $\{b = 0, a \geq 0\}$ do not belong to $C(4,0)$: indeed, the segment $r(f)$ coincides with the intersection set of relation curves of all polynomials of entire pencil $\{b = 0\}$. Therefore for such $f$ $\tilde{\pi}^{-1}(f, \dots, f) \cap (F_3 \setminus F_2) = \emptyset$.

**Proposition 9.** *For any cell $X$ other than $C$ and $O$, the set $\{X\}$ is the space of a locally trivial fibre bundle over our cell, whose fiber over $f$ is an open cone over*

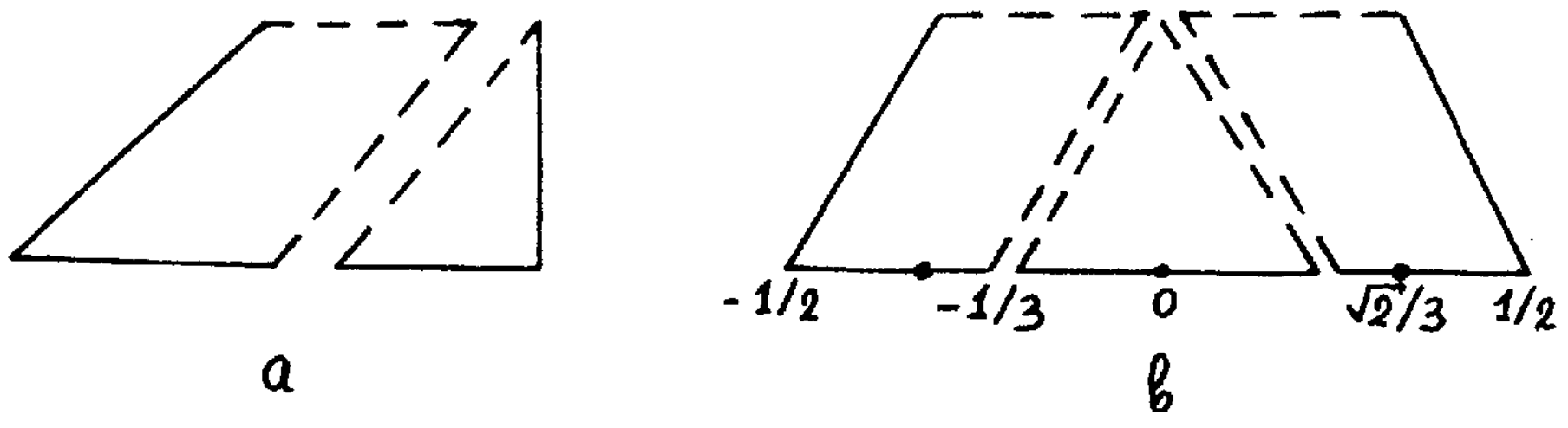

Fig. 4: The set of lines through $f$, $f \in B$ or $D$, satisfying three elementary conditions

*the set* $\tilde{\pi}^{-1}(f, \dots, f) \cap F_2$.

For any cell this will follow from the forthcoming explicit description of this set.

### 4.4.2. Cells $A$ and $A'$.

**Lemma 9.** *If a polynomial* $f = t^4 + at^2 + bt$ *has only one real critical point, then* $\tilde{\pi}^{-1}(f, \dots, f) \cap (F_2 \setminus F_{1.5}) = \emptyset$.

This follows immediately from the Corollary of Lemma 7. □

**Corollary.** *For any* $f$ *from the cell* $A$ *or* $A'$, *the set* $\tilde{\pi}^{-1}(f, \dots, f) \cap (F_3 \setminus F_2)$ *is homeomorphic to the open cone over the space* $S^{*2}$ *(see § 4.3.2) where* $S$ *is the relation curve of* $f$ *(in particular,* $S$ *is homeomorphic to a line segment).* □

**Corollary.** *The closed homology group of the part* $\{A\}$ *or* $\{A'\}$ *of* $F_3 \setminus F_2$ *swept out by the sets considered in the previous corollary is trivial in all dimensions.*

This follows immediately from the Künneth formula and Lemma 8. □

### 4.4.3. The cells $B$ and $B'$.

**Proposition 10.** *For any* $f$ *from the cells* $B$ *and* $B'$, *the set of affine lines* $(f, g) \subset P_4$ *through* $f$, *such that all polynomials of the line satisfy exactly three common elementary conditions, consists of two components and is homeomorphic to the figure 4a) (where the dashed curve denotes the part of the boundary not belonging to the set).*

The proof of this proposition and of Propositions 11, 12 below is given in § 5.

**Corollary.** *The closed homology groups* $\bar{H}_*(\{B\})$ *and* $\bar{H}_*(\{B'\})$ *are trivial.*

*Proof.* By Proposition 9 and Künneth formula, it is sufficient to prove that for any $f$ from our cell the closed homology group of any set $\tilde{\pi}^{-1}(f, \dots, f) \cap F_2$ is trivial. By Proposition 10, this homology coincides with that of $\tilde{\pi}^{-1}(f, \dots, f) \cap F_{1.5}$ (while the similar homology of $\tilde{\pi}^{-1}(f, \dots, f) \cap (F_2 \setminus F_{1.5})$ is trivial). By construction of $F_{1.5}$, the latter group coincides with that of $S^{*2}$ where $S$ is the relation curve of $f$, i.e. the disjoint union of two line segments. This implies the Corollary. $\Box$

#### 4.4.4. The cells $E$ and $E'$.

**Proposition 11.** *For any $f$ from the cells $E$ and $E'$, the set of affine lines $(f, g)$ in $P_4$ passing through $f$ and such that all polynomials of the line satisfy exactly three common elementary conditions is homeomorphic to the half-open interval.*

**Corollary.** *The closed homology groups $\bar{H}_*(\{E\})$ and $\bar{H}_*(\{E'\})$ are trivial.*

The proof is the same as for the Corollary of Proposition 10 (in this case the relation curve of $f$ is the disjoint union of a line segment and a point.) $\Box$

#### 4.4.5. The cell $D$.

**Proposition 12.** *For any $f$ from the cell $D$, the set of affine lines $(f, g)$ in $P_4$ passing through $f$ and such that all polynomials of the line satisfy exactly three common elementary conditions consists of three connected components and is homeomorphic to the figure 4b) (where the dashed segment denotes the part of the boundary which does not belong to the set and the double dashed segments inside the figure should be distracted from it).*

**Corollary.** *The closed homology group $\bar{H}_*(\{D\})$ is trivial.*

The proof is the same as for the Corollary of Proposition 10 (for $f \in D$ the curve $r(f)$ is homeomorphic to the closed cross and thus $r(f)^{*2} \sim \{a\ point\}$.) $\Box$

Theorem 3 follows immediately from Corollaries of Lemma 9 and of Propositions 10, 11 and 12 and terminates the proof of Theorem 1.

## 5. Proof of Propositions 10, 11 and 12

Let $f = t^4 + at^2 + bt$ be a point of the cell $B$ or $B'$. Let $x_1 < x_2 < x_3$ be three roots of the polynomial $f' \equiv 4t^3 + 2at + b$; obviously $x_1 + x_2 + x_3 = 0$. It is sufficient to consider only one of the cells $B, B'$, say $B$. In this case $x_1$ is negative, $x_2$ and $x_3$ are positive, and the topological shape of the relation curve $r(f)$ is as shown in Fig. 5a), so that $x_2$ and $x_3$ are connected by the "finite" branch of it. The line $\{t + s = 0\}$ is asymptotical for $r(f)$, in particular it does not intersect it at finite points and separates two branches of it.

By the Corollary of Lemma 7, if $\#(r(g) \cap r(f)) = 3$, then the topological disposition of these curves is as in Fig. 2a) or 2b) (or is some degeneration of this figure corresponding to the collision of some endpoint(s) of $r(g)$ with these of $r(f)$).

**Notation.** For any $f \in B$, denote by $\heartsuit \equiv \heartsuit(f)$ (respectively, by $\diamondsuit$) the set of

$$g = t^3 + \alpha t^2 + \beta t \tag{5.1}$$

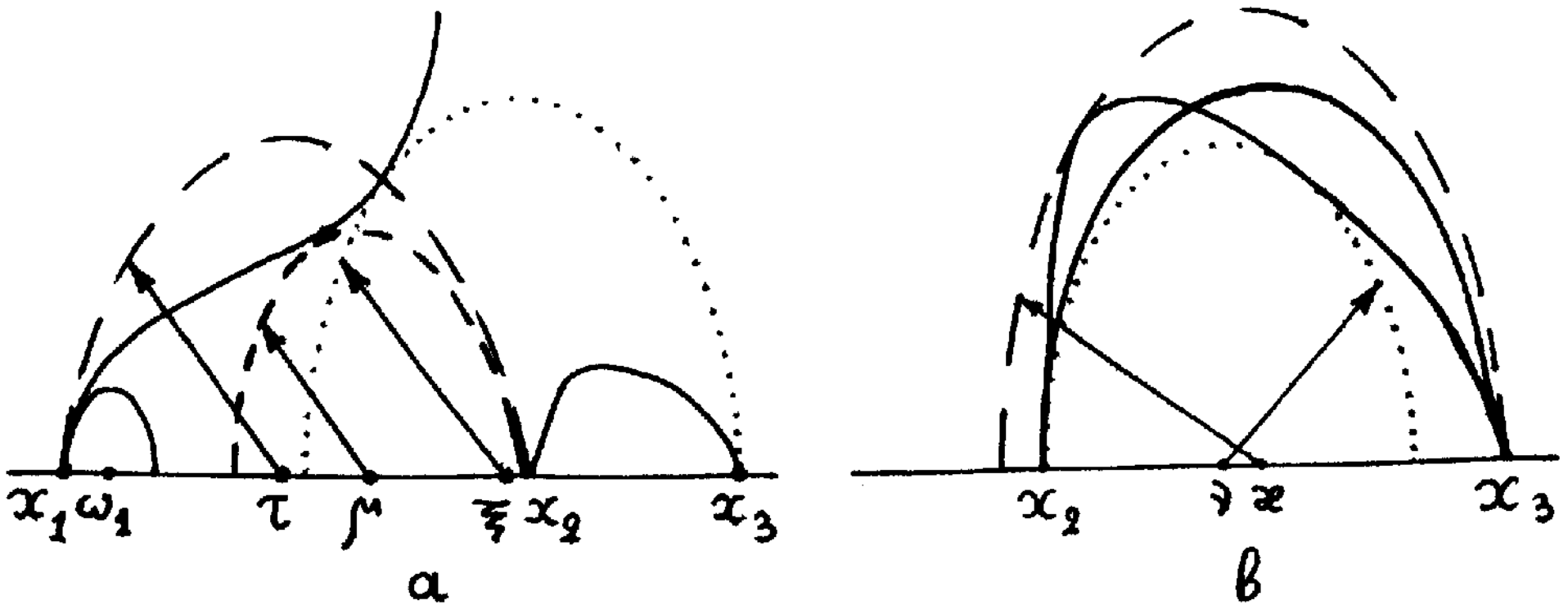

Fig. 5: Relation curve for $f$ of type $B$ and critical positions of $r(g)$

such that the intersection $r(f) \cap r(g)$ consists of three points, two of which belong to the "infinite" branch of $r(f)$ (respectively, all three points belong to the "finite" branch).

**Theorem 4.** *For any $f \in B$ the set $\heartsuit$ (respectively, $\diamondsuit$) is homeomorphic to the right-hand (respectively, left-hand) part of Fig. 4a).*

The proof of this theorem occupies the next subsections 5.1–5.3; Proposition 10 is a direct corollary of it.

## 5.1. Families $L^\nabla$ and focal points

**Notation.** For any $\nabla \in \mathbf{R}^1$ denote by $L^\nabla$ the line in $P_3$ consisting of polynomials (5.1) with $\alpha = -3\nabla$, or, which is the same, with center of $r(g)$ at $(\nabla, \nabla)$.

The curves $r(g)$, $g \in L^\nabla$, are the level sets of the function

$$F_\nabla \equiv (t^2 + ts + s^2) - 3\nabla(t + s). \tag{5.2}$$

**Lemma 10.** *The restriction of a quadratic form in $\mathbf{R}^2$ on a cubic curve has at most 9 isolated critical points.*

This follows immediately from the Bezout's theorem. □

**Corollary.** *The restriction of any function $F_\nabla$ on $r(f)$, $f \in B$, has at most three critical points in the interior points of $r(f) \subset \Psi(4)$.*

Indeed, the curve (4.7) is symmetric with respect to the involution $\{t \leftrightarrow s\}$, and three fixed points $x_i$ of it are the critical points of any $F_\nabla$. □

**Lemma 11.** *There are values $\omega_1 > x_1$, $\omega_2 > x_2$ and $\omega_3 < x_3$ such that if $\nabla$ lies on the same side of $\omega_i$ as $x_i$ (respectively, $\omega_i$ separates $\nabla$ and $x_i$) then the* *restriction of $F_\nabla$ on $r(f)$ has at $x_i$ a strong (quadratic) local minimum (respectively, maximum). Namely,*

$$\omega_1 = \frac{6x_1^2 - a}{12x_1}, \quad \omega_2 = \frac{6x_2^2 - a}{12x_2}, \quad \omega_3 = \frac{6x_3^2 - a}{12x_3}. \tag{5.3}$$

Indeed, $(\omega_i, \omega_i)$ are the focal points of $(r(f), (x_i, x_i))$ in the sense of [7] in the Euclidean metric $t^2+ts+s^2$. The precise values (5.3) follow from direct calculations. □

## 5.2. The structure of the set ♡.

For any $\nabla \in \mathbf{R}^1$ we consider the metamorphoses of the curves $r(g)$ along the line $L^\nabla$. For $\beta = 3\nabla^2$ this curve consists of single point $(\nabla, \nabla)$, but when the parameter $\beta$ decreases this curve blows around the same center. Its disposition with respect to $r(f)$ can undergo the following four metamorphoses essential for us:

$M_i$ ($i = 1, 2, 3$): the intersection point of $r(g)$ and the diagonal $\{t = s\}$ passes through the intersection point $(x_i, x_i)$ of $r(f)$ with the diagonal;

$M_4$ (optional) the curve $r(g)$ touches the "infinite" (i.e., lying in the domain $\{t + s < 0\}$) branch of $r(f)$ in some interior point of it (and intersects it in the forthcoming instants corresponding to greater ellipses or smaller values of $\beta$.)

Let us denote the values of $\beta$ at which these metamorphoses happen by $\beta_1, \dots, \beta_4$ respectively; of course they are functions of $\nabla$ (or of $\alpha \equiv -3\nabla$) and of the coefficients $a, b$ of $f$. $g$ belongs to the set ♡ if and only if the following conditions are satisfied:

a) the metamorphose $M_4$ actually happens in the family $L^\nabla$ at some instant $\beta_4(\alpha) > \beta$;

b) $\beta \geq \beta_1(\alpha)$; $\beta_2(\alpha) \geq \beta$; $\beta \geq \beta_3(\alpha)$.

For a given $\alpha$, the set of values $\beta$ for which the polynomial (5.1) satisfies the conditions a), b) can be empty, or consist of unique point (if $\beta_4(\alpha) > \beta_2(\alpha) = \beta_1(\alpha) > \beta_3(\alpha)$), be a closed segment (if $\beta_4(\alpha) > \beta_2(\alpha) > \max(\beta_1(\alpha), \beta_3(\alpha))$) or a half-open interval (if $\beta_2(\alpha) > \beta_4(\alpha) > \max(\beta_1(\alpha), \beta_3(\alpha))$).

Now we investigate the values of $\alpha$ for which these conditions are satisfied.

**Lemma 12.** *The metamorphose $M_4$ actually happens at some instant $\beta_4(\alpha) > \beta_1(\alpha)$ if and only if $\nabla > \omega_1$.* □

### 5.2.1. Breakpoints

Now for any given $f \in B$ we define three important points in the line $\mathbf{R}^1$ of the values of $\nabla$, see Fig. 5a).

Set $\tau = (x_1 + x_2)/2$, so that $(\tau, \tau)$ is the center of the ellipse of the form (4.6) intersecting the diagonal at the points $\{x_1, x_1\}$ and $\{x_2, x_2\}$.

The family of ellipses (4.6) intersecting the diagonal at two points, the bigger of which (i.e. the closest to $(+\infty, +\infty)$) is $(x_3, x_3)$, is obviously a one-parametric curve in the plane $\{\alpha, \beta\}$; all ellipses of this family lie inside one another. Consider the smallest ellipse of this family intersecting the "infinite" branch of $r(f)$, and denote by $\xi$ the corresponding value $-\alpha/3$, so that $(\xi,\xi)$ is the center of this ellipse. Finally, the point $\mu$ is defined in similar way through the family of ellipses (4.6) containing the point $(x_2,x_2)$.

**Lemma 13.** $\tau > \omega_1$.

Indeed, by Lemma 11 (and since $x_1 < 0$) this assertion is equivalent to the inequality $6x_1x_2 + a < 0$. But for $f \in B$ both $a$ and $x_1x_2$ are negative. □

**Corollary 1.** $\mu > \tau$.

Indeed, by definition $\mu \geq \tau$. But if $\mu = \tau$ then the ellipse $r(g)$ through $(x_2,x_2)$ with center at $(\mu,\mu)$ intersects the diagonal also at the point $(x_1,x_1)$. By Lemma 13 it bounds a continuum of points of the "infinite" branch, which contradicts the definition of $\mu$. □

**Corollary 2.** $\xi > \mu$.

Indeed, if $\xi \leq \mu$ then the ellipse (4.6) through $(x_2,x_2)$ with center at $(\mu,\mu)$ (participating in the definition of $\mu$) lies strictly inside the ellipse through $(x_3,x_3)$ with center $(\xi,\xi)$. Since the former has nonempty intersection with the "infinite" branch of $r(f)$, the latter bounds a continuum of points of this branch, which contradicts the definition of $\xi$. □

Finally, we obtain that the breakpoints are ordered as follows: $\omega_1 < \tau < \mu < \xi$. The line $L^\nabla$ intersects the domain $\heartsuit$ by the set homeomorphic to:
- a point if $\nabla = \tau$,
- a line segment if $\nabla \in (\tau,\mu)$,
- a half-open interval if $\nabla \in [\mu,\xi)$,
- empty set if $\nabla \notin [\tau,\xi)$.

This implies the assertion of Theorem 4 about $\heartsuit$.

**Remark.** In this proof we considered the fibres of the map $\heartsuit \to \mathbf{R}^1$ sending $g$ to the center of symmetry of $g'$. Using another map $\heartsuit \to \mathbf{R}^1$ sending $g$ to its greater critical point we can establish a homeomorphism of $\heartsuit$ to the direct product of a line segment and a half-open interval. Indeed, by Corollary of Lemma 7 this critical point should belong to the segment $[x_2,x_3]$ and by Corollary 1 of Lemma 13 for any point $x$ of this segment the fibre of our map over $x$ is a half-open interval.

## 5.3. The structure of the set $\diamondsuit$.

**Lemma 14.** *If the line $L^\nabla$ has a non-empty intersection with the set $\diamondsuit$, then $\nabla \leq (x_2+x_3)/2$.*

Indeed, this follows from the Corollary of Lemma 7. □

**Lemma 15.** $(x_2+x_3)/2 < \omega_2$ *and* $(x_2+x_3)/2 < \omega_3$.

Indeed, by (5.3) (and since $x_3 > 0 < x_2$) both these inequalities are equivalent to the inequality $6x_2x_3 + a < 0$. By the Vieta formulae

$$x_1 + x_2 + x_3 = 0, \quad x_1x_2 + x_1x_3 + x_2x_3 = a/2, \tag{5.4}$$

this is equivalent to $(x_3 - x_2)^2 > 0$. $\Box$

By the Roll theorem, if $r(g)$ has three intersection points with the "finite" branch of $r(f)$, $g \in L^\nabla$, then the restriction of the function $F_\nabla$ on this branch has at least two critical points in its interior part.

**Lemma 16.** *Suppose that $\nabla \leq (x_2+x_3)/2$. Then the restriction of $F_\nabla$ has two critical points inside the "finite" branch if and only if $\nabla > \sqrt[3]{b}/2$.*

*Proof.* The condition $dF_\nabla \| df$ can be easily calculated and is equivalent in the interior of $\Psi(4)$ to the equation $(t^2 + 4ts + s^2) + 2\alpha(t + s) - a = 0$. Together with the equation of $r(f)$ this gives the following system in the symmetric functions $p, q$:

$$\begin{cases} p(p^2 - 2q) + ap + b = 0 \\ p^2 + 2q + 2\alpha p - a = 0 \\ p^2 - 4q \geq 0, \end{cases} \tag{5.5}$$

cf. (4.10). Two first equations in it give the following equation on $p$:

$$p^3 + \alpha p^2 + b/2 = 0. \tag{5.6}$$

If $\nabla = (x_2+x_3)/2$, then $F_\nabla$ has two critical points inside the "finite" component (by Lemma 15 and since $F_\nabla$ takes equal values at the endpoints of this component); for $\nabla$ close to $-\infty$ $F_\nabla$ obviously has no such critical points. When $\nabla$ moves from $(x_2+x_3)/2$ to $-\infty$ the number of critical points can decrease only when the discriminant of the polynomial (5.6) becomes zero or when the last inequality in (5.5) becomes an equality. The latter is impossible: indeed, it corresponds to the situation when $t = s$, i.e. some interior critical point of $F_\nabla|_{r(f)}$ comes to an endpoint of the component, which contradicts Lemma 15. And the former discriminant condition is obviously equivalent to $\nabla = \sqrt[3]{b}/2$. $\Box$

**Lemma 17.** $x_2 < \sqrt[3]{b}/2 < (x_2 + x_3)/2$.

Indeed, both these inequalities follow from the Vieta formulae (5.4) and $b/4 = -x_1x_2x_3$ and the fact that $x_3 > x_2$ and $-x_1 > 2x_2$. (The right-hand inequality follows also from the proof of the previous lemma.) $\Box$

So, for any point $\nabla \in (\sqrt[3]{b}/2, (x_2 + x_3)/2]$ the restriction of $F_\nabla$ on the finite branch has four critical points: the minima at the points $x_2$ and $m(\nabla)$ and the maxima at $x_3$ and $M(\nabla)$. Denote by $\beta^1(\nabla), \ldots, \beta^4(\nabla)$ minus the critical values of these points respectively, so that all ellipses $r(g)$, $g = t^3 - 3\nabla t^2 + \beta^i(\nabla)t$ touch the finite branch of $r(f)$. The intersection of our line $L^\nabla$ with $\Diamond$ consists of the polynomials (5.1) with $\alpha = -3\nabla$ and $\beta$ satisfying both conditions $\beta^1(\nabla) \geq \beta \geq$ $\beta^3(\nabla)$ and $\beta^2(\nabla) > \beta > \beta^4(\nabla)$. Thus to investigate the set $\diamondsuit$ we need only to count the order of points $\beta^i(\nabla)$ for different $\nabla \in (\sqrt[3]{b}/2, (x_2+x_3)/2]$.

For any such $\nabla$ the critical values of both maxima points are obviously greater than these of minima, hence $\beta^1(\nabla)$ and $\beta^2(\nabla)$ are greater than $\beta^3(\nabla)$ and $\beta^4(\nabla)$.

**Lemma 18.** *For any $f \in B$, there is a point $\kappa \in (\sqrt[3]{b}/2, (x_2+x_3)/2)$ such that $\beta^3(\kappa) = \beta^4(\kappa)$, $\beta^3(\nabla) > \beta^4(\nabla)$ if $\nabla > \kappa$ and $\beta^3(\nabla) < \beta^4(\nabla)$ if $\nabla < \kappa$.*

Indeed, consider the one-parametric family of polynomials (5.1) such that one of their critical points coincides with $x_3$ and the other is smaller or equal than $x_2$. This family intersects any line $L^\nabla$, $\nabla \leq (x_2+x_3)/2$, at exactly one point, thus $\nabla$ is a natural coordinate in it. Let $g_0$ be the initial point of the family (which belongs to $L^\nabla$, $\nabla = (x_2+x_3)/2$), then the ellipse $r(g_0)$ intersects the finite branch of $r(f)$ at three points, two of which are $(x_2, x_2)$ and $(x_3, x_3)$, see Fig. 5b). When $\nabla$ decreases, the ellipses $r(g)$ of the family grow monotonically, in particular starting with some instant $\kappa$ they contain this branch of $r(f)$ inside. At this instant the corresponding ellipse touches the branch at two points (one of which is $(x_3, x_3)$) and the critical values $-\beta^3(\kappa), -\beta^4(\kappa)$ of the corresponding function $F_\kappa$ coincide. By the Roll theorem, $F_\kappa$ has one more critical point inside our branch, hence by Lemma 16 $\kappa > \sqrt[3]{b}/2$. Obviously this value $\kappa$ is the desired one. $\square$

**Lemma 19.** *For any $f \in B$, there is a point $\nu \in (\sqrt[3]{b}/2, (x_2+x_3)/2)$ such that $\beta^1(\nu) = \beta^2(\nu)$, $\beta^1(\nabla) > \beta^2(\nabla)$ if $\nabla < \nu$ and $\beta^1(\nabla) < \beta^2(\nabla)$ if $\nabla > \nu$.*

The proof repeats that of Lemma 18: we consider the family of polynomials $g$ with one critical point equal to $x_2$ and the other smaller or equal to $x_3$ (so that the initial function $g_0$ is the same as in the previous proof, see Fig. 5b)). $\square$

The topological structure of the set $\diamondsuit$ follows from Lemmas 14–19. More precisely, the intersection of the line $L^\nabla$ with this set is a point if $\nabla = (x_2+x_3)/2$, a closed segment if $\nabla \in (\max(\kappa,\nu), (x_2+x_3)/2)$, a half-open interval if $\nabla \in (\min(\kappa,\nu), \max(\kappa,\nu)]$, an interval if $\nabla \in (\sqrt[3]{b}/2, \max(\kappa,\nu)]$, and is empty if $\nabla \notin (\sqrt[3]{b}/2, (x_2+x_3)/2]$. $\square$

## 5.4. Proof of Proposition 11.

The set investigated in Proposition 11 is the limit of sets $\heartsuit = \heartsuit(f)$ when $f \in B$ tends to the cell $E$ and hence the segment $[x_2, x_3]$ contracts to a point. The structure of it follows from the Corollary 1 of Lemma 13. $\square$

## 5.5. Proof of Proposition 12.

The group of scalings (4.9) acts effectively on the cell $D$ and acts freely on the corresponding part of $\tilde{\rho}\Sigma$. Therefore we can consider unique polynomial $f = t^4 - 2t^2$, so that $r(f)$ consists of the segment

$$t + s = 0 \qquad (5.7)$$

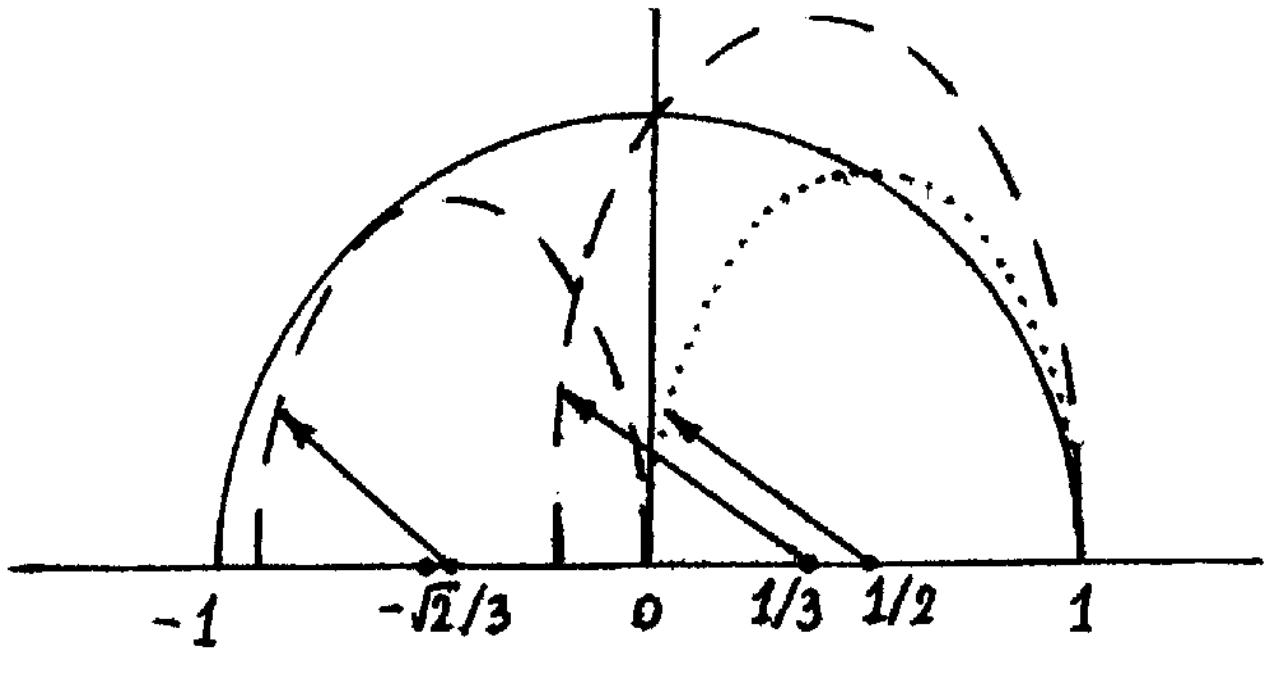

Fig. 6: Intersections with the symmetric relation curve

and half-circle

$$t^2 + s^2 = 2. \tag{5.8}$$

**Lemma 20.** *The set of such values of $\nabla$ that for some $g \in L^\nabla$ the relation curve $r(g)$ intersects the half-circle (5.8) twice is the interval* $(-2/3, 2/3)$.

This is the limit version of Lemma 11 (more precisely, of its part concerning the points $x_1$ and $x_3$) when $f \in B$ tends to our polynomial $t^4 - 2t^2 \in D$. □

If $\nabla \in (1/2, 2/3] \cup [-2/3, -1/2)$, then the ellipses of the family $L^\nabla$ pass through the endpoint $(1,1)$ or $(-1,-1)$ of the half-circle (5.8) earlier than they meet the component (5.7), hence such families also do not intersect the investigated set.

For $\nabla = \pm 1/2$ the intersection of the family $L^\nabla$ with this set consists of unique polynomial $t^3 \mp (3/2)t^2$, see Fig. 6, and for $|\nabla|$ slightly smaller than $1/2$ it consists of a segment: indeed, the metamorphoses of the dispositions of ellipses $r(g)$ of such families are ordered as follows: first (for the greatest value of $\beta$) they touch the half-circle (5.8) at its interior point, then meet for the first time the segment (5.7), and then pass through the endpoint $(1,1)$ or $(-1,-1)$ of the half-circle; the segment of the family connecting the second and the third occasions is the desired one.

The situation changes if the first and the second metamorphoses happen simultaneously or even change their order: in this case the boundary point at which the family enters the interesting set does not belong to this set. It is easy to calculate that this is the case if $|\nabla| \leq \sqrt{2}/3$.

Moreover, any family $L^\nabla$ has an exceptional point, the ellipse $r(g)$ passing through the crucial point $(-1,1)$ of $r(f)$, which intersects $r(f)$ by a smaller number of points than the neighboring ellipses $r(g')$. If $|\nabla| \leq 1/3, \nabla \neq 0$, then this exceptional point lies inside the interval of interesting points of the family and thus should be distracted from it; the union of such points is depicted by the double dashed segments in Fig. 4b). □